\documentclass[11pt]{amsart}
\usepackage{geometry}   
\usepackage{setspace}
\usepackage{graphicx}
\usepackage{amssymb, amsthm}
\usepackage{epstopdf}
\usepackage{tikz}
\usetikzlibrary{matrix,arrows}
\usepackage{array}
\usepackage{txfonts}
\usepackage{stmaryrd}

  \newcommand{\Z}{\ensuremath{{\mathbb{Z}}}}

  \newcommand{\G}{\Gamma}
  \newcommand{\iso}{\cong}
\newcommand{\GG}{\ensuremath{{G_\G}}}
\newcommand{\tGG}{\ensuremath{{G_{\tilde\G}}}}
\newcommand{\FG}{\ensuremath{{F_\G}}}

\newcommand{\CATZ}{{\rm CAT}(0)}

\newtheorem{theorem}{Theorem}[section]
    \newtheorem{corollary}[theorem]{Corollary}
   \newtheorem{lemma}[theorem]{Lemma}

    \theoremstyle{definition}
\newtheorem{definition}[theorem]{Definition}
\newtheorem{remark}[theorem]{Remark}

\newtheorem{examples}[theorem]{Examples}

\title{The automorphism group of a graph product with no SIL}
\author{Ruth Charney  \and Kim Ruane \and Nathaniel Stambaugh  \and Anna Vijayan}
\thanks {R. Charney was partially supported by NSF grant DMS 0705396.}  

\begin{document}

\begin{abstract}
\noindent
We study the automorphisms of graph products of cyclic groups, a class of groups that includes all right-angled Coxeter and right-angled Artin groups.  We show that the group of automorphism generated by partial conjugations is itself a graph product of cyclic groups providing its defining graph does not contain any separating intersection of links (SIL). In the case that all the cyclic groups are finite, this implies that the automorphism group is virtually CAT(0); it has a finite index subgroup which acts geometrically on a right-angled building. 
 \end{abstract}
 
 \maketitle

\section{Introduction}

Classically, a collection of groups can be combined using free products or direct products.  More generally, a graph product of groups is a class of groups which interpolates between these. 
Let $\Gamma$ be a finite simplicial graph with vertex set $V$ and let $\{G_v\}_{v\in V}$ be a family of groups.  Then the graph product $G_{\Gamma}$ is the quotient of the free product of the groups $G_v$ obtained by adding commutator relations between $G_v$ and $G_w$ whenever $v,w$ are adjacent in $\Gamma$.  A discrete graph $\Gamma$ gives the free product of the $G_v$ and a complete graph gives  the direct product.  
Graph products encompass several important classes of groups.  In particular, one obtains the class of right-angled Coxeter groups by requiring each $G_v$ be isomorphic to $\mathbb Z\slash 2\mathbb Z$ and the class of right-angled Artin groups when each $G_v$ is isomorphic to $\mathbb Z$.   In this paper we require only that the vertex groups be finitely generated abelian groups. Any such graph product  is isomorphic to a graph product of cyclic groups, hence we can restrict our attention to the latter.

The automorphism groups of right-angled Coxeter and right-angled Artin groups have been studied extensively in the literature (see, for example, \cite{Ti}, \cite{Mu}, \cite{GuKa}, \cite{La}, \cite{ChCrVo}, \cite{ChVo}, \cite{BuChVo}, \cite{Day}).  Automorphisms of more general graph products were considered by Laurence in his thesis \cite{La-th}.    Building on work of Servatius \cite{Se},   Laurence describes a finite generating set for $Aut(G_{\Gamma})$ in the case when all vertex groups have the same order, either infinite or a fixed prime $p$.  More recently, in \cite{GuPiRu},  Gutierrez, Piggott, and Ruane   begin a unified treatment of the automorphism group of a general graph product of cyclic groups, and in \cite{CoGu}, Corredor and Gutierrez extend Laurence's generating set to all such graph products.

The automorphism group of $\GG$ is generated by four types of automorphisms:  graph symmetries, vertex isomorphisms,  transvections, and partial conjugations.  The first two types generate a finite subgroup.  Transvections, which map $v \mapsto vw$ (or $v \mapsto wv$) for a pair of vertices $v,w$, are familiar to those who work with free group automorphisms. Indeed, the automorphism group of a free group is entirely generated by these transvections.  For graph products, however, the transvections are more restricted (in some cases excluded entirely) and the partial conjugations play an essential role. 
 
 Partial conjugations are defined as follows. Given a vertex $v\in V$, let $lk(v)$ denote the full subgraph of $\Gamma$ spanned by the vertices  adjacent to $v$ and $st(v)$ the subgraph spanned by $v$ and $lk(v)$.  For each connected component $C$ of $\Gamma\setminus st(v)$, the {\it partial conjugation} $\pi_{v,C}$ conjugates all of the generators in $C$ by $v$ and leaves all other generators fixed. 

The subgroup of $Aut(G_{\Gamma}$) generated by the partial conjugations is denoted $Aut^{pc}(G_{\Gamma})$ and will be our main object of study.  In the case where all vertex groups have finite order, this subgroup has finite index in the full automorphism group $Aut(G_{\Gamma})$.  This is also true for some graph products with infinite vertex groups, namely those for which the structure of $\G$ does not permit any transvections (eg., if $\G$ has no circuits of length less than four and no valence one vertices).

% It is easy to show that some partial conjugations commute with each other.  Our main theorem (Theorem \ref{main}) is states that under appropriate conditions on the graph $\G$, $Aut^{pc}(G_{\Gamma})$ is itself a graph product of cylcic groups with the defining graph determined by these commutator relations.

A simplicial graph $\Gamma$ has a {\it Separating Intersection of Links (SIL)} if for some pair $v,w$ with $d_{\Gamma}(v,w) \ge 2$, there is a component of $\Gamma\setminus (lk(v)\cap lk(w))$ which contains neither $v$ nor $w$.  Our main theorem, Theorem \ref{main}, states that if $\G$ has no SILs, then  $Aut^{pc}(G_{\Gamma})$ is itself a graph product of cyclic groups.

%\begin{theorem}
%If $G_\Gamma$ is a graph product of cyclic groups, and $\Gamma$ has no SILS, then $Aut^{pc}(G_{\Gamma})$ is also a graph product of cyclic groups of the same orders.
%\end{theorem}

To prove this, we consider the graph $\tilde\Gamma$ whose vertices are in one-to-one correspondence with the partial conjugations $\pi_{v,C}$ of $G_{\Gamma}$.  Two vertices of $\tilde\Gamma$ are connected by an edge if the two partial conjugations commute, thus we have a graph product of cyclic groups $G_{\tilde\Gamma}$.  In the case where $\Gamma$ has no SILS, we prove that $Aut^{pc}(G_{\Gamma})$ is isomorphic to $G_{\tilde\Gamma}$.   The main technical point is to characterize exactly when two partial conjugations commute.  Under the no SILS assumption, we give a simple characterization of when $\pi_{v,C}$ and $\pi_{w,D}$ commute in terms of the relative position of $C$ and $D$. This is the content of Lemma \ref{comrel}.

Our main theorem has some interesting geometric implications. 
Recall that a \emph{CAT(0) metric space} is a proper, complete metric space in which each geodesic triangle is ``at least as thin'' as the Euclidean triangle with the same  side lengths.  We say that a finitely generated group $G$ is a \emph{CAT(0) group} if $G$ acts properly, cocompactly by isometries on a CAT(0) metric space (such an action is said to be {\it geometric}).
A group $G$ is {\it virtually CAT(0)} if some finite index subgroup of $G$ is CAT(0). Note that extending a geometric action from a finite index subgroup to the full group is highly non-trivial.  It is unknown if virtually CAT(0) groups are CAT(0).  

In Section \ref{prelim}, we show that  any graph product of cyclic groups $\GG$ acts on a right-angled building.  Right-angled buildings are always  CAT(0) by a theorem of Davis \cite{Da2}.  If the vertex groups are all finite, then the associated building is locally finite, its automorphism group is a locally compact group, and the graph product sits as a uniform lattice in this group.   In particular,  $\GG$ is a CAT(0) group.  (For a discussion of right-angled buildings and their lattices see \cite{Th} and \cite{BaTh}.)   If some vertex group is infinite cyclic, then the action is no longer proper.  However, 
if all the vertex groups are infinite (the right-angled Artin group case), then there is a different CAT(0) cube complex, the Salvetti complex, which can be used to get a CAT(0) structure on $\GG$.
 
Our main theorem implies that under the no SILs hypothesis, $Aut^{pc}(G_{\Gamma})$ is itself a graph product hence also acts on a right-angled building.  Moreover, we show that this action extends to the larger group generated by partial conjugations, graph symmetries and vertex isomorphisms.    If all of the vertex groups are finite, these generate a finite index subgroup of $Aut(\GG)$ and we conclude that $Aut(\GG)$ is virtually CAT(0).  (This last statement also follows from \cite{GuPiRu} where they show that under these hypotheses, the inner automorphism group has finite index in $Aut(\GG)$. Our construction gives a CAT(0) action of a much larger subgroup, sometimes encompassing the entire automorphism group.)  

One would like to know whether, in general, these actions can be extended to the full automorphism group $Aut(G_{\Gamma})$, that is, whether the action can be extended to include transvections.  This would almost certainly require a different geometric construction as transvections do not behave well with respect to the geometry of the cube complexes given here.  Even in the case of right-angled Coxeter and Artin groups, it is unknown whether the full automorphism groups are $\CATZ$.

%%%%%%%%%%%%%%%%%%%%%%%%%%%%%%%%%%%

\section{Graph products and associated geometries}\label{prelim}

In this section we discuss graph products and their associated geometries.  We begin with a definition of a graph product.

\begin{definition}\label{gr-product}
Let $\G$ be a finite, simplicial graph with vertex set $V$, together with a labeling of each vertex by a group $G_v$.  Let $\FG$ denote the free product of all the vertex groups $G_v, v \in V$. Then the \emph{graph product} $\GG$ is the quotient group of $\FG$ obtained by adding commutator relations between $G_v$ and $G_w$ whenever $v,w$ are connected by an edge in $\G$.  
\end{definition}

In this paper we investigate graph products of cyclic groups, that is, graph products for which all of the vertex groups $G_v$ are \emph{cyclic}.  More generally, if all of the vertex groups of a graph product $\GG$ are finitely generated abelian groups, then $\GG$ is naturally isomorphic to the graph product obtained by replacing each vertex in $\G$ by a complete graph with vertices labelled by the (indecomposable) cyclic summands of $G_v$.  Thus, our results apply more generally to this class of groups.

Gutierrez and Piggott \cite{GuPi}, generalizing work of Laurence \cite{La-th},  have shown that for any graph product of indecomposable cyclic groups, the graph $\G$ and the vertex groups $G_v$ are uniquely determined by the isomorphism class of $\GG$. Thus, when referring to the graph product $\GG$, we may assume that this data has been specified.  

\medskip
\emph{For the remainder of the paper, we assume that all vertex groups are cyclic.}
\medskip

\begin{examples}  If all of the vertex groups are cyclic of order 2, then we obtain the right-angled Coxeter groups.  If all of the vertex groups are infinite cyclic, then we obtain the right-angled Artin groups.
\end{examples}

Given $g \in \GG$, a \emph{reduced word} for $g$ is a minimal length word $g_1g_2\dots g_k$ in $\FG$  (with each $g_i$ belonging to some vertex group) representing $g$.   Any word representing $g$ can be reduced by a process of ``shuffling" (i.e., interchanging commuting elements) and combining adjacent elements from the same vertex group. Any two reduced words representing $g$ differ only by shuffling \cite{Gr-th}.

For any subset $T$ of the vertex set $V$, let $G_T$ denote the graph product associated to the full subgraph of $\G$ spanned by $T$.  The natural map from $G_T$ into $\GG$ splits, hence $G_T$ is isomorphic to its image and we make no distinction between them.  By convention, we set $G_{\emptyset}={1}$. 

To a graph product $\GG$, we associate a cubical complex $X_\G$ as follows.  Define two sets, partially  ordered by inclusion, 
\begin{align*}
\mathcal S_\G &= \{ G_T \mid T \subseteq V, ~\text{$G_T$ is abelian}\}\\
&\iso \{ T \mid T \subseteq V, ~\text{$T$ spans a complete subgraph of $\G$}\}\\
G\mathcal S_\G &= \{ gG_T \mid g \in \GG, ~T \subseteq V, ~\text{$G_T$ is abelian}\}.
\end{align*}
Let $X_\G$ be the geometric realization of the poset $G\mathcal S_\G$ and let $K \subset X_\G$ be the geometric realization of $\mathcal S_\G$. 
Left multiplication of $\GG$ on this poset induces an action of $\GG$ on $X_\G$.  A fundamental domain for this action is $K$, and hence the action is cocompact.  The stabilizer of the vertex $gG_T$ is conjugate to $G_T$ which is finite if and only if all the vertex 
groups in $T$ are finite.  Thus, the action of $\GG$ on $X_\G$ is proper if and only if $\G$ is a graph of \emph{finite} cyclic groups.

The complexes $X_\G$ are interesting in their own right.  As we will now show, they have the structure of right-angled buildings.  These buildings are based on a construction of Davis  \cite{Da3}, \cite{Da2}.  
In the case of a right-angled Coxeter group, $X_G$ is the well-known Davis complex.  For a right-angled Artin group, $X_\G$ is known as the Deligne complex (or in the terminology of \cite{ChDa}, the ``modified" Delinge complex).  We follow \cite{Da} and \cite{BaTh} for basic definitions.  

First recall that a \emph{chamber system over a set $S$} is a set $\Phi$ of \emph{chambers} together with a family of equivalence relations on $\Phi$ indexed by $S$.  For $s \in S$, we say two chambers are \emph{$s$-adjacent} if they are $s$-equivalent, but not equal.  For a word $w=s_1\dots s_k$, $s_i \in S$, a \emph{gallery} of type $w$ is a sequence of chambers $\phi_0,\phi_1, \dots ,\phi_k$ such that $\phi_{i-1}$ is $s_i$-adjacent to $\phi_i$.

Now suppose $S$ is the generating set of a right-angled Coxeter group $W$.  A \emph{$W$-valued distance function} on $\Phi$ is a function $d : \Phi \times \Phi \to W$ such that, given a reduced word $s_1\dots s_k$ representing $w \in W$, there exists a gallery of type $s_1\dots s_k$ from $\phi$ to $\phi'$ if and only if $d(\phi,\phi')= w$.

\begin{definition}\label{bldg}
Let $W=W_\G$ be a right-angled Coxeter group with generating set $S$. Then a \emph{right-angled building of type $W$} is a chamber system $\Phi$ over $S$ such that
\begin{enumerate}
\item for all $s \in S$, every $s$-equivalence class contains at least two chambers,
\item there exists a $W$-valued distance function $d : \Phi \times \Phi \to W$.
\end{enumerate}
\end{definition}

Let $\GG$ be a graph product of cyclic groups.  Denote by $W_\G$ the right-angled Coxeter group obtained by replacing each vertex group $G_v$ by $W_v = \Z/2\Z$.  Define a set-theoretic map (not a homomorphism) $\gamma : \GG \to W_\G$ as follows.  For $g \in \GG$, represent $g$ by a reduced word $g=g_1 \dots g_k$, with $g_i \in G_{v_i}$, and set $\gamma(g)=s_1 s_2 \dots s_k$ where $s_i$ is the generator of $W_{v_i}$.  This is well-defined since any two reduced words for $g$ are related by commutator relations which also hold in $W_\G$.  Moreover, $s_1 \dots s_k$ is also reduced since no shuffling of $g_1 \dots g_k$ (and hence of $s_1 \dots s_k$) allows two elements of the same vertex group to be combined.

\begin{theorem}\label{XGbldg}  For any graph product $\GG$ of cyclic groups, $X_\G$ is a right-angled building of type $W_\G$.
\end{theorem}

\begin{proof}  We take $\Phi$ to be the set of translates of $K$ in $X_\G$ and we say that two chambers $gK, hK$ are $s_i$-equivalent if $g^{-1}h \in G_{v_i}$.   Then every $s_i$-equivalence class contains $q$ elements where $q = |G_{v_i}|$.  

Define $d : \Phi \times \Phi \to W_\G$ by $d(gK, hK)= \gamma(g^{-1}h)$. Then for a reduced word 
$w=s_1s_2\dots s_k$, there exists a gallery of type $w$ from $gK$ to $hK$   if and only if
$g^{-1}h= g_1 \dots g_k$ for some $g_i \in G_{v_i}$, or equivalently,  $d(gK,hK)=w$.  
\end{proof}

These buildings and their automorphism groups are studied by Barnhill, Thomas, Haglund, and Paulin \cite{Th},  \cite{BaTh}, \cite{Ha}, \cite{HaPa}.  If the vertex groups are all finite, then the (full) automorphism group of the building is a locally compact topological group and $\GG$ is a uniform lattice in this group.

Although $X_\G$ was defined as a simplicial complex, it has a natural cubical structure whose cubes correspond to  ``intervals".  For a pair of subsets $T_1 \subseteq T_2$ in $\mathcal S_\G$, the interval $[G_{T_1},G_{T_2}]$ is the subcomplex of $K$ spanned by the vertices $G_T$, $T_1 \subseteq T \subseteq T_2$.  It is combinatorially a cube of dimension $|T_2 - T_1|$.  The translates of these intervals give a cubical structure on all of $X_\G$. 

The fundamental chamber $K$ is independent of the orders of the vertex groups.  Thus, it is isometric to the fundamental chamber in the Davis complex for $W_\G$.  It was shown by Davis in \cite{Da2} that any such right-angled building is CAT(0) with respect to the cubical metric described above.  (This can also be proved directly for $X_\G$ using the link condition for cubical complexes.)

The action of $\GG$ takes intervals to intervals, hence preserves the cubical metric and the quotient by $\GG$ is just the fundamental chamber $K$.   Thus $\GG$  acts faithfully (the stabilizer of $G_{\emptyset}$ is trivial), cocompactly by isometries on $X_\G$.   As noted above, however, the action is proper if and only if every vertex group is finite.

\begin{corollary} For all graph products of cyclic groups, the cubical metric on $X_\G$ is CAT(0). If the vertex groups are all finite, then $\GG$ is a CAT(0) group.
\end{corollary}

\begin{remark}\label{action}  For use later in the paper, we remark that this action can be extended to a slightly larger group.  Let $\Sigma_\G$ be the (finite) group of automorphisms of $\GG$ generated by symmetries of the graph $\G$ (which permute the generators of $\GG$) and automorphisms of a single vertex group.  This group acts on the poset $G \mathcal S_\G$ in the obvious way,
$\sigma \cdot gG_T = \sigma(gG_T)$,  and hence it acts on $X_\G$.  Combining this with the $\GG$-action gives an action of the semi-direct product $\GG \rtimes \Sigma_\G$ on $X_\G$.  This action is again proper, cocompact, isometric, and faithful.
 \end{remark}

%%%%%%%%%%%%%%%%%%%%%%%%%%%%%%%%%%%%

\section{Automorphism groups and separating intersections of links}

In this section we introduce the no SILs condition on $\G$ and study automorphism groups of graph products of cyclic groups $\GG$ satisfying this condition.  

Servatius \cite{Se} and Laurence \cite{La}, \cite{La-th} described a finite generating set for $Aut(\GG)$ for certain classes of graph products, such as right-angled Artin groups.  This result has recently been extended to all graph products of cyclic groups by Corredor and Gutierrez in \cite{CoGu}.  We now describe this generating set.

In order to simplify notation, we will think of the vertex $v$ as the generator of the cyclic group $G_v$, so that the vertex set $V$ generates $\GG$.
Denote the order of $v$ (and hence of $G_v$) by $|v|$.  Associated to a vertex $v$ in $\G$ are two subgraphs:  the link of $v$, $lk(v)$, is the full subgraph spanned by the vertices adjacent to $v$ and the star of $v$, $st(v)$, is the subgraph spanned by $v$ and $lk(v)$.  

\begin {theorem}[\cite{La}, \cite{CoGu}]  \label{generate} 
If $\GG$ is a graph product of cyclic groups, then $Aut(\GG)$ is generated by automorphisms of the following types:

\begin{enumerate}
 \item   \emph{Symmetries}: induced by symmetries of $\G$, permute the generators

 \item \emph{Vertex isomorphisms}: automorphisms of a single vertex group $G_v$
            
  \item  \emph{Partial conjugations}: conjugate all of the generators in one connected component  $C$ of $\G \smallsetminus st(v)$ by $v$
  
   \item   \emph{Transvections}:  map  $v\mapsto vw^k$ or $v \mapsto w^kv$ where one of the following holds
   \begin{enumerate}
      \item $|v|=\infty$, $k=1$, and $lk(v) \subseteq st(w)$, or
    \item $|v|=p^i$, $|w|=p^j$,   $k=\max\{1, p^{j-i}\}$, and $st(v) \subseteq st(w)$.
    \end{enumerate}  
\end{enumerate}
\end{theorem}

We are interested primarily in the partial conjugations.  Denote by $\pi_{v,C}$  the partial conjugation by $v$ of the component C, and let  $Aut^{pc}(W_\Gamma)$ denote  the group generated by all partial conjugations. 
 
It follows from Lemma 2.8 of \cite{GuPiRu} that  when the vertex groups are all finite, $Aut^{pc}(W_\Gamma)$ has finite index in the full automorphism group $Aut(\GG)$.  This is also the case when there are no permissible transvections (for example when $\G$ has no cycles of length less than 5 and no vertices of valence less than 2).  

The interaction between two partial conjugations $\pi_{v,C}$ and $\pi_{w,D}$ depends on the relative position of the components $C$ and $D$.
A crucial role will be played by the following.

\begin{definition} 
A simplicial graph $\Gamma$ has a \emph{Separating Intersection of Links (SIL)} if for some pair $(v,w)$, with $d_\Gamma(v,w)\geq 2$, there is a component of $\Gamma \smallsetminus (lk(v) \cap lk(w))$ which contains neither $v$ nor $w$.
\end{definition}

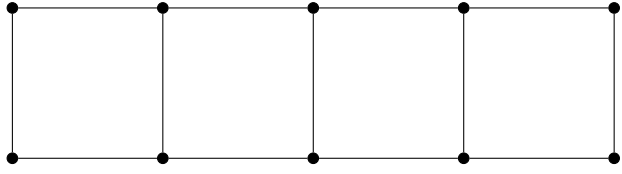
\begin{figure}
\begin{center}
\begin{tikzpicture}[scale=2] 
\tikzstyle{every node}=[circle, draw, fill=black!100,                       inner sep=0pt, minimum width=4pt]    
\path (0,0cm) node [label=0:$ $] (a1) {};  
\path (1,0cm) node [label=0:$ $] (a2) {};  
\path (2,0cm) node [label=0:$ $] (a3) {};  
\path (3,0cm) node [label=0:$ $] (a4) {};  
\path (4,0cm) node [label=0:$ $] (a5) {};  
\path (0,1cm) node [label=0:$ $] (b1) {};  
\path (1,1cm) node [label=0:$ $] (b2) {};  
\path (2,1cm) node [label=0:$ $] (b3) {};  
\path (3,1cm) node [label=0:$ $] (b4) {};  
\path (4,1cm) node [label=0:$ $] (b5) {};    
\draw (a1) -- (a2) -- (a3) -- (a4) -- (a5) -- (b5) -- (b4) -- (b3) -- (b2) -- (b1) -- (a1); \draw (a2) -- (b2); 
\draw (a3) -- (b3); 
\draw (a4) -- (b4); 
\end{tikzpicture} 
\end{center} 
\caption{A graph with separating stars but no SILs}  
\label{noSIL} 
\end{figure}

\begin{remark}\label{disconnected} The no SILs condition is most interesting for connected graphs. For if $\G$ has more than two connected components, then it necessarily has a SIL.  If it has two components $\G_1$ and $\G_2$, one of which is not a complete graph, then it also has a SIL since if $v$ and $w$ are vertices in $\G_1$ with $d(v,w) \geq 2$, then $\G_2 \subset \Gamma \smallsetminus (lk(v) \cap lk(w))$ is a component containing neither $v$ nor $w$.  Thus a graph with no SILs is either connected or it is the disjoint union of two complete graphs.  
\end{remark}

In the case where $\G$ has no SILs, we will prove that $Aut^{pc}(\GG)$ is itself a graph product of cyclic groups $\tGG$ where the vertex set of $\tilde\Gamma$ corresponds to the set of partial conjugations. The edges of $\tilde\Gamma$ will correspond to the partial conjugations that commute, and are prescribed by the following lemma.

\begin{lemma} \label{comrel}
Suppose $\Gamma$ is a connected simplicial graph which does not contain any SILs and let $v$ and $w$ be vertices of $\Gamma$. Suppose $d(v,w)\geq 2$, and let $C_0$ be the component of $\Gamma \smallsetminus st(v)$ containing $w$, and $D_0$ be the component of $\Gamma \smallsetminus st(w)$ containing $v$. Then
\begin{enumerate}
\item Every component of $\Gamma \smallsetminus st(v)$, except $C_0$, lies entirely in $D_0$, and every component of $\Gamma \smallsetminus st(w)$, except $D_0$, lies entirely in $C_0$.
\item The partial conjugations $\pi_{v,C}$ and $\pi_{w,D}$ commute unless $C=C_0$ and $D=D_0$.
\end{enumerate}
\end{lemma}

\begin{proof}
(1) Let $C$ be a component of $\Gamma \smallsetminus st(v)$. If $C$ contains any vertex of $lk(w)$, then it also contains $w$, so $C = C_0$. If not, then $C \cap st(w) = \emptyset$, so $C$ lies completely in some component $D$ of $\Gamma \smallsetminus st(w)$. We claim that $D=D_0$. 

Let $\bar{C}$ denote the graph generated by $C$ and the vertices adjacent to $C$. Clearly $\bar{C} \smallsetminus C \subset lk(v)$. On the other hand, $\bar{C} \smallsetminus C \not\subset lk(v) \cap lk(w)$, since this would imply that $C$ was a component of $\Gamma \smallsetminus lk(v) \cap lk(w)$ which did not contain $v$ or $w$. It follows that $C$ and $v$ are adjacent to a vertex which is not in the link of $w$. Hence $C$ and $v$ are in the same component of $\Gamma \smallsetminus st(w)$, i.e., $D=D_0$ as claimed.

(2) First we note that $\pi_{v,C_0}$ and $\pi_{w,D_0}$ do not commute by direct computation.
\begin{align*}
\pi_{v,C_0} \circ \pi_{w,D_0}(v) &= \pi_{v,C_0} (wvw^{-1}) = vwvw^{-1}v^{-1}\\
\pi_{w,D_0} \circ \pi_{v,C_0}(v) &= \pi_{w,D_0} (v) = wvw^{-1}
\end{align*}
Next consider the case where $C \neq C_0$ and $D \neq D_0$. By (1) $C \cap D = \emptyset$, and we do another direct computation.
$$\pi_{v,C} \circ \pi_{w,D}(x) = \pi_{w,D} \circ \pi_{v,C}(x) = \begin{cases} vxv^{-1} & x \in C \\ wxw^{-1} & x \in D \\ x & x \not\in (C\cup D) \end{cases}$$
Now suppose $C \neq C_0$ and $D = D_0$. Then by (1), we know that $C\subset D$, $v \in D$, and $w \not\in C$. We can once again check this by direct computation.
$$\pi_{v,C} \circ \pi_{w,D}(x) = \pi_{w,D} \circ \pi_{v,C}(x) = \begin{cases} wvxv^{-1}w^{-1} & x \in C \\ wxw^{-1} & x \in D\smallsetminus C \\ x & x \not\in D \end{cases}$$
The remaining case where $C = C_0$ and $D \neq D_0$ is similar.
\end{proof}

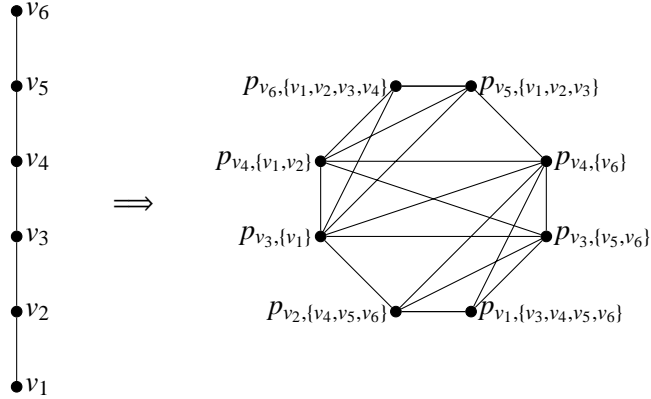
\begin{figure}
\begin{center}
\begin{tabular}{m{1cm}m{1cm}m{6cm}}
\begin{tikzpicture}[scale=1]
\tikzstyle{every node}=[circle, draw, fill=black!100,
                       inner sep=0pt, minimum width=4pt]
   \path (0,0cm) node [label=0:$v_1$] (a1) {}; 
	\path (0,1cm) node [label=0:$v_2$] (a2) {}; 
	\path (0,2cm) node [label=0:$v_3$] (a3) {}; 
	\path (0,3cm) node [label=0:$v_4$] (a4) {}; 
	\path (0,4cm) node [label=0:$v_5$] (a5) {}; 
	\path (0,5cm) node [label=0:$v_6$] (a6) {};

	\draw (a1) -- (a2) -- (a3) -- (a4) -- (a5) -- (a6);
\end{tikzpicture}
&$\Longrightarrow$&
\begin{tikzpicture}[scale=1]
\tikzstyle{every node}=[circle, draw, fill=black!100,
                       inner sep=0pt, minimum width=4pt]
   \path (.5,1cm) node [label=0:$p_{v_1,\{v_3,v_4,v_5,v_6\}}$] (v1) {}; 
	\path (-.5,1cm) node [label=180:$p_{v_2,\{v_4,v_5,v_6\}}$] (v2) {}; 
	\path (-1.5,2cm) node [label=180:$p_{v_3,\{v_1\}}$] (v3a) {}; 
	\path (1.5,2cm) node [label=0:$p_{v_3,\{v_5, v_6\}}$] (v3b) {}; 
	\path (-1.5,3cm) node [label=180:$p_{v_4,\{v_1, v_2\}}$] (v4a) {}; 
	\path (1.5,3cm) node [label=0:$p_{v_4,\{v_6\}}$] (v4b) {}; 
	\path (.5,4cm) node [label=0:$p_{v_5,\{v_1,v_2,v_3\}}$] (v5) {}; 
	\path (-.5,4cm) node [label=180:$p_{v_6,\{v_1,v_2,v_3,v_4\}}$] (v6) {}; 
	
	\draw (v3b) -- (v1) -- (v2) -- (v3a) -- (v4a) -- (v5) -- (v6) -- (v3a);
	\draw (v2)--(v4b) -- (v1) -- (v2) -- (v3b) -- (v4b) -- (v5) -- (v6) -- (v4a);
	\draw (v3a) -- (v4b) -- (v4a) -- (v3b) -- (v3a) --(v5);
	
%	\draw[thin,dashed,red] (v4a) -- (v1) -- (v5) -- (v2) -- (v6) -- (v1) -- (v3a) -- (v5) -- (v3b) -- (v6);
%	\draw[thin,dashed,red] (v4b) -- (v2);
\end{tikzpicture}
\end{tabular}
\end{center}
\caption{The graphs $\G$ and $\tilde\G$}\label{graphs}
\end{figure}

We now construct the graph $\tilde{\Gamma}$.  The vertices of $\tilde\Gamma$ are in one-to-one correspondence with the partial conjugations $\pi_{v,C}$, and are denoted by $\tilde V=\{p_{v,C}\}$.  Any two vertices $p_{v,C}$ and $p_{w,D}$ are connected by an edge unless $d(v,w)\geq 2$, $v \in D$, and $w \in C$.  We assign to the vertex  $p_{v,C}$ the cyclic group of  order $|v|$.  An example is shown in Figure \ref{graphs}.

By the lemma, there is a homomorphism $\phi: \tGG \rightarrow Aut^{pc}(\GG)$ which takes $p_{v,C} \mapsto \pi_{v,C}$.  This homomorphism is clearly surjective; our goal is to prove that if $\G$ contains no SILs, then $\phi$ is an isomorphism.  To do this, we will pass to the outer automorphism group. 

The outer automorphism group of $\GG$ is the quotient of $Aut(\GG)$ by the subgroup $Inn(\GG)$ of inner automorphsims of $\GG$. The inner automorphism by a vertex $v$ is the product of all the partial conjugations by $v$, hence this subgroup lies in $Aut^{pc}(\GG)$ and we can define $Out^{pc}(\GG)$ accordingly.  We would like to define a corresponding quotient for $\tGG$.

The inner automorphism group is isomorphic to the group modulo its center.  In the case of 
$\GG$, the center is generated by the vertices (if any) which are connected to every other vertex in $\G$.  Let $\Delta$ be the (possibly empty) graph generated by these vertices and let $\G_0 = \G \smallsetminus \Delta$. Then $\GG$ decomposes as the direct product of $G_{\Delta}$ and $G_{\G_0}$, so $Inn(\GG)$ is isomorphic to $G_{\G_0}$. 

We denote by $p_v$ the product $p_v = \prod p_{v,C}$ over all components $C$ of 
$\G \smallsetminus st(v)$, so that $\phi(p_v)$ is the inner automorphism by $v$. 

\begin{lemma} \label{normal}
The correspondence $v \mapsto p_v$ induces a homomorphism $\tilde{f}:G_{\G_0} \rightarrow \tGG$ and the image of  $\tilde{f}$ is a normal subgroup of $\tGG$. 
\end{lemma}
\begin{proof}
The first statement follows by definition of $\tilde\G$ since if $d(v,w) \leq 1$ in $\G$, then $p_{v,C}$ commutes with $p_{w,D}$ for all $C,D$, so commuting relations are preserved and the order of $p_v$ is $|v|$. 

To prove that the image is normal, we will show that for any $v \in \G_0$ and for any generator of $p_{w,D}$ of $\tGG$, the following equation holds.
\begin{equation*}
p_{w,D}p_vp^{-1}_{w,D} = \begin{cases} p_v & \text{ if } v \not\in D \\
			          p_wp_vp_w^{-1} & \text{ if } v \in D \end{cases}
\end{equation*}
Note that in either case, $p_{w,D}p_vp^{-1}_{w,D}$ is in $\tilde f(G_{\G_0})$. 

{\it Case 1:} $v \not\in D$. Then by Lemma \ref{comrel} $p_{w,D}$ commutes with $p_{v,C}$ for every $C$. 

{\it Case 2:} $v \in D$. Consider the expression $p_w p_v p_w^{-1}$. Then for each connected component $D^{\prime}$ of $\Gamma \smallsetminus st(w)$ with $D \neq D^{\prime}$, the partial conjugation $p_{w,D^{\prime}}$ commutes with $p_v$ by Lemma \ref{comrel}. Simplifying the expression, we get the desired result.
\end{proof}

In light of the lemma, we can now form the quotient group, $Q = \tGG/ \tilde f(G_{\G_0})$.
If $f$ denotes the inclusion of the inner automorphisms into $Aut^{pc}(\GG)$, then the diagram below clearly commutes.

\begin{center}
\begin{tikzpicture}[description/.style={fill=white,inner sep=2pt}]
\matrix (m) [matrix of math nodes, row sep=3em,
column sep=2em, text height=1.5ex, text depth=0.25ex]
{ G_{\Gamma_0} &  & \tGG \\
Inn(\GG) &  & Aut^{pc}(\GG) \\ };
\path[->,font=\scriptsize]
(m-1-1) edge node[auto] {$ \tilde{f} $} (m-1-3)
(m-1-1) edge node[auto] {$ \cong $} (m-2-1)
(m-2-1) edge node[auto] {$ f $} (m-2-3)
(m-1-3) edge node[auto] {$ \phi $} (m-2-3);
\end{tikzpicture}
\end{center}

It follows that $\tilde f$ is injective and that $\phi$ induces a map on the quotient groups, so we have 
a commutative diagram  of exact sequences,

\begin{center}
\begin{tikzpicture}[description/.style={fill=white,inner sep=2pt}]
\matrix (m) [matrix of math nodes, row sep=3em,
column sep=2em, text height=1.5ex, text depth=0.25ex]
{ 1 & & G_{\Gamma_0} &  & G_{\tilde{\Gamma}} && Q && 1\\
1 & & Inn(G_\Gamma) &  & Aut^{pc}(G_\Gamma) && Out^{pc}(G_\Gamma) && 1 \\ };
\path[->,font=\scriptsize]
(m-1-1) edge node[auto] {} (m-1-3)
(m-1-3) edge node[auto] {$ \tilde{f} $} (m-1-5)
(m-1-5) edge node[auto] {$ \tilde{g} $} (m-1-7)
(m-1-7) edge node[auto] {} (m-1-9)
(m-2-1) edge node[auto] {} (m-2-3)
(m-1-3) edge node[auto] {$ \cong $} (m-2-3)
(m-2-3) edge node[auto] {$ f $} (m-2-5)
(m-2-5) edge node[auto] {$ g $} (m-2-7)
(m-2-7) edge node[auto] {} (m-2-9)
(m-1-5) edge node[auto] {$ \phi $} (m-2-5)
(m-1-7) edge node[auto] {$ \psi $} (m-2-7);
\end{tikzpicture}
\end{center}

We are now ready to state and prove our main result. 

\begin{theorem}\label{main} Let $\GG$ be a graph product of cyclic groups whose defining graph $\G$ contains no SILs.  Then the map $\phi: \tGG \to Aut^{pc}(\GG)$ is an isomorphism.  In particular, $Aut^{pc}(\GG)$ is a graph product of cyclic groups of the same order(s) as the vertex groups of $\GG$.
\end{theorem}

\begin{proof}
In light of the exact sequence above, it suffices to prove that the map $\psi$ on the quotient groups is an isomorphism. The theorem will then follow from the 5-lemma.  Since $\phi$ and $g$ are both surjective, $\psi \circ \tilde{g}  = g \circ \phi$ is surjective, and so $\psi$ is as well.

Suppose $\Gamma$ is connected. We first argue that $Q$ is abelian. Take two generators $p_{v,C}$ and $p_{w,D}  \in \tGG$ which {\it do not} commute. By Lemma \ref{comrel},  we know that $p_{v,C}$ {\it does} commute with $p_{w,D^{\prime}}$ for every connected component $D^{\prime}$ of $\Gamma \smallsetminus st(w)$ with $D^{\prime} \neq D$, hence it commutes with the product  $p^{\prime} = \prod_{D^{\prime}\neq D} p_{w,D^{\prime}}$. But $p^{\prime}$ and $p_{w,D}$ represent inverse elements in $Q$, so the images of $p_{v,C}$ and $p_{w,D}$ commute in $Q$.  

Since $Q$ is abelian, for any $\bar x \in  Q$, we can write $\bar x$ as a product of generators $p_{v,C}$ in which all occurrences of a vertex $v$ appear together.  That is, we can choose a representative $x \in \tGG$ of the form
$$x = \prod_v \prod_C p_{v,C}^{k_{v,C}}.$$
Now suppose $\bar x \in ker(\psi)$.  By the commutivity of the relevant diagram, $\phi(x)
= \prod_v \prod_C \pi_{v,C}^{k_{v,C}}$ lies in  $Inn(W_\Gamma)$.  A product of this form is an inner automorphism if and only if, for fixed $v$, the power  $k_{v,C}$ is the same for every component $C$ of $\G \smallsetminus st(v)$. This means that  $x$ has the form $\prod_v p_v^{k_v}$, which is clearly in $G_{\G_0}$, thus $\bar x$ is trivial in $Q$.  We conclude that the kernel of $\psi$ is trivial, so $\psi$ is an isomorphism.

It remains to consider the case when $\Gamma$ is not connected. By Remark \ref{disconnected}, this occurs only when $\Gamma$ is the disjoint union of two complete graphs.  In this case, every partial conjugation is an inner automorphism, so $\G = \tilde\G $, and since $\GG$ has trivial center, $Aut^{pc}(\GG)= Inn(\GG) \cong \GG$. 
\end{proof}

As noted in Section 2, any graph product of finitely generated abelian groups is isomorphic to a graph product of cyclic groups obtained by ``blowing up" a vertex $v$ with group $G_v$ into a complete graph with vertices labeled by the indecomposable cyclic summands of $G_v$.  Applying Theorem \ref{main}
to this new graph product, we see that $Aut^{pc}(G_\Gamma)$ is again a graph product of cyclic groups.  Moreover, this graph product is just the blow-up of  $\tGG$ (where $\tilde\G$ is defined as above, but the vertex groups $G_v$ are no longer cyclic).  Thus we may restate the theorem as follows.

\begin{theorem}\label{general}
Let $\GG$ be a graph product of finitely generated abelian groups whose defining graph $\G$ contains no SILs.  Then the map $\phi: \tGG \to Aut^{pc}(\GG)$ is an isomorphism.  In particular, 
 $Aut^{pc}(G_\Gamma)$ is also a graph product of finitely generated abelian groups. 
 \end{theorem}

We remark that the proof of the theorem above also gives an independent proof of the following result of \cite{GuPiRu}.

\begin{corollary}[\cite{GuPiRu}]  Assume $\G$ is connected. Then $Out^{pc}(G_\Gamma)$ is abelian if and only if $\G$ contains no SILs.
\end{corollary}

\begin{proof}  In the proof of the main theorem we showed that if $\G$ has no SILs, then $Q$, and hence  $Out^{pc}(\GG)$, is abelian.  If $\G$ has a SIL, that is a component $C$ of 
$\Gamma \smallsetminus (lk(v) \cap lk(w))$ containing neither $v$ nor $w$ where $d(v,w)\geq 2$, then it is straightforward to check that the commutator $[\pi_{v,C}, \pi_{w,C}]$ is not an inner automorphism.
Hence $Out^{pc}(\GG)$ is not abelian.
\end{proof}

%%%%%%%%%%%%%%%%%%%%%%%%%%%%%%%%%%

\section {Geometric implications}

Recall that a group $G$ is a CAT(0) group if it acts geometrically (i.e., properly, cocompactly by isometries) on a complete CAT(0) space.
It is an open question whether automorphism groups of graph products, even in the Coxeter group case, are CAT(0) groups.  Theorem \ref{main} gives some partial answers.

In the case where all vertex groups are finite, $Aut^{pc}(\GG)$ has finite index in $Aut(\GG)$ so by Theorems \ref{XGbldg} and \ref{main} we obtain 

\begin{corollary}\label{action1}
Let $\GG$ be a graph product of finite cyclic groups whose defining graph has no SILs. Then the automorphism group $Aut(\GG)$ is virtually CAT(0).  More precisely,  there is a faithful, geometric action of  $Aut^{pc}(\GG)$ on the right-angled building $X_{\tilde\G}$.
\end{corollary}

It would be nice to extend this action to the whole automorphism group.  Recall from Theorem \ref{generate} that $Aut(\GG)$ is generated by four types of automorphisms: symmetries, vertex isomorphisms, partial conjugations, and transvections.
 Letting $\Sigma_\G$ denote the group generated by symmetries and vertex isomorphisms, the subgroup of $Aut(\GG)$ generated by the first three types of automorphisms 
is a semi-direct product, $Aut^{pc}(\GG) \rtimes \Sigma_\G$.
 We can easily extend the action of $Aut^{pc}(\GG)$ on $X_{\tilde\G}$ to this larger group.

\begin{corollary}\label{action2}
 Let $\GG$ be a graph product of finite cyclic groups whose defining graph has no SILs. Then the action of $Aut^{pc}(\GG)$ on  $X_{\tilde\G}$ extends to a faithful, geometric action of
$Aut^{pc}(\GG) \rtimes \Sigma_\G$.
\end{corollary}

\begin{proof}
 By Remark \ref{action}, the action of $\tGG$ on $X_{\tilde\G}$   extends to a faithful, geometric action of the semi-direct product $\tGG \rtimes \Sigma_{\tilde\G}$.  The group $\Sigma_\G$ embeds naturally in $\Sigma_{\tilde\G}$ (an isomorphism of $G_v$ goes to the product of the corresponding isomorphisms of $G_{p_{v,C}}$ for all components $C$).  Combining this embedding with the isomorphism $\phi^{-1}$, we get an inclusion $Aut^{pc}(\GG) \rtimes \Sigma_\G \hookrightarrow \tGG \rtimes \Sigma_{\tilde\G}$, and hence an induced action on $X_{\tilde\G}$.
\end{proof}

If some of the vertex groups are infinite cyclic, then the action of $\tGG$ on $X_{\tilde\G}$ is not proper.  However, if \emph{all} of the vertex groups are infinite, then $\GG$ and $\tGG$ are right-angled Artin groups  and we can use a different geometric construction, the Salvetti complex (see \cite{Ch}), to get an action on a CAT(0) space.

\begin{corollary}\label{action3}
Let $\Gamma$ be a simplicial graph with no SILs, and suppose $\GG$ is a right-angled Artin group. Then the subgroup of $Aut(\GG)$  generated by partial conjugations, inversions and graph symmetries acts faithfully and geometrically on a CAT(0) cube complex, the Salveti complex of $\GG$. 
\end{corollary}

\begin{proof}  It is easy to show that the action of $A_{\tilde\G}$ on its Salvetti complex extends to an action of  $A_{\tilde\G}  \rtimes \Sigma_{\tilde\G}$.  The proof then proceeds as above.  
\end{proof}

We close by remarking that some graph products $\GG$ of cyclic groups do not permit transvections in which case the subgroup in Corollaries \ref{action2} and \ref{action3}  constitutes the entire automorphism group.  This is the case, for example, if $\G$ has no triangles and no vertices of valence less than two, or if every pair of adjacent vertex groups have relatively prime order.  For those that do permit transvections,  the action described above does not extend in any natural way to an isometric action of the transvections.  In this case, proving that the full automorphism group is CAT(0) will almost certainly require a different space.

\end{document}